\documentclass[10pt,a4paper,twoside]{amsart}

\usepackage{amsfonts, amssymb, amsmath, amsthm}
\usepackage{mathrsfs} 
\usepackage{latexsym}
\usepackage{enumerate}
\usepackage{verbatim}

\usepackage{url}
 
\usepackage[colorlinks=true]{hyperref}
\hypersetup{citecolor=blue, linkcolor=blue}

\newtheorem{thm}{Theorem}[section]
\newtheorem{lem}[thm]{Lemma}

\def\goes{\mathrel{\rightarrow}}
\def\Ocal{\mathcal{O}}
\def\1{1\!\!\!1}

\title{A Moebius sum%
}
\author{Olivier Ramar\'e}

\begin{document}
 \address[O. Ramar\'e]{CNRS/ Institut de Math\'ematiques de Marseille, Aix 
 Marseille Universit\'e, U.M.R. 7373, Site Sud, Campus de Luminy, Case 907, 
 13288 
 Marseille Cedex 9, France.}
 \email{olivier.ramare@univ-amu.fr}

\date{\sl August 15th, 2023}


\maketitle

\begin{abstract}
  We provide numerical bounds for $\Sigma(X)=\sum_{\substack{d_1,d_2\le
        X}}\frac{\mu(d_1)\mu(d_2)}{[d_1,d_2]}$. We show in particular
    that $0\le \Sigma(X)\le 17/25$ for every $X\ge2$.
\end{abstract}


\section{Introduction and results}

In several problems of analytic number theory appears the sum
\begin{equation}
  \sum_{\substack{d_1,d_2\le X}}\frac{\mu(d_1)\mu(d_2)}{[d_1,d_2]}.
\end{equation}
It has been shown in \cite{Dress-Iwaniec-Tenenbaum*83} by F.~Dress,
G.~Tenenbaum and H.~Iwaniec, that this sum
converges to a constant, and H.~Helfgott computes precisely
its first four decimals in Proposition 6.30 of
\cite{Helfgott*30}. This implies in particular
 that this constant is $\le 0.4408$.
Our aim in this note is to provide a uniform and explicit
upper bound for this sum.
\begin{thm}
  \label{Main}
  We have
  \begin{equation*}
    0\le \sum_{\substack{d_1,d_2\le
        X}}\frac{\mu(d_1)\mu(d_2)}{[d_1,d_2]}
    \le
    \begin{cases}
      17/25&\text{when $X\ge 2$},\\
      0.574&\text{when $X\ge 10^{9}$},\\
      0.536&\text{when $X\ge 3\cdot 10^{10}$},\\
      0.504&\text{when $X\ge 2.4\cdot 10^{12}$}.
    \end{cases}
  \end{equation*}
  We have $17/25=0.68$.
\end{thm}
\noindent
The bounds for $X$ have been chosen with an application in mind.
The lower estimate comes from the interpretation of this constant as
follows. By the Parseval equality, we have that, when $\sigma>1$:
\begin{align*}
  \lim_{T\goes\infty}
  \frac{1}{2i\pi T}\int_{-\infty}^\infty
  \biggl|\zeta(s)\sum_{n\le X}\frac{\mu(d)}{d^s}\biggr|^2ds
  &=\sum_{n\ge 1}
  \biggl(\sum_{\substack{d|n\\ d\le X}}\mu(d)\biggr)^2/n^{2\sigma-1}
  \\&=\zeta(2\sigma-1)\sum_{\substack{d_1,d_2\le X}}
  \frac{\mu(d_1)\mu(d_2)}{[d_1,d_2]^{2\sigma-1}}
\end{align*}
and a continuity argument ends the proof.
\subsection*{Notation}
Our notation is classical: we use
\begin{equation}
  \label{eq:4}
  m_q(y)=\sum_{\substack{d\le y\\ (d,q)=1}}\frac{\mu(d)}{d},
  \quad
  \text{and}
  \quad
  m(y)=m_1(y).
\end{equation}
Furthermore, $f=\Ocal^*(g)$ means that $|f|\le g$.

\section{On the Moebius function}

\begin{lem}
  \label{m1}
  We have $|m(x)|\le \sqrt{2/x}$ for $0<x\le 10^{14}$ and $|m_2(x)|\le
  \sqrt{3/x}$ for $0<x\le 10^{12}$.
\end{lem}
\noindent
This is \cite[Lemma 5.10]{Helfgott*30} and
\cite[Eq. (5.79)]{Helfgott*30}.
\begin{lem}
  \label{m2}
  We have $|m(x)|\le 0.0144/\log x$ for $x \ge 463\,421$ and $|m_2(x)|\le
  0.0296/\log x$ for $x\ge 5379$.
\end{lem}
\noindent
This follows from \cite[Theorem 1.2]{Ramare*12-5} and from \cite[Lemma 5.17]{Helfgott*30}

\begin{lem}
  \label{m3}
  With $\xi=1-1/(12\log 10)$, $y>1$ and any $t\in(0,y]$,  we have
  \begin{equation*}
    |m(t)|\le \sqrt{2/t}
    +0.0144\cdot\1_{y\ge 10^{12}}\frac{y^{1-\xi}}{\log y}\frac{1}{t^{1-\xi}}
  \end{equation*}
  as well as
  \begin{equation*}
    |m_2(t)|\le \sqrt{3/t}
    +0.0296\cdot\1_{y\ge 10^{12}}\frac{y^{1-\xi}}{\log y}\frac{1}{t^{1-\xi}}
  \end{equation*}
\end{lem}
\noindent
This is a trivial modification of \cite[Lemma 5.12]{Helfgott*30}: we
degraded the value of $\xi$ in the first estimate to get a more
uniform result, added the condition $\1_{t\ge 10^{12}}$ which was
obvious in the proof and transposed almost verbatim the argument to $m_2$.
\begin{lem}
  \label{m4}
  With $\xi=1-1/(12\log 10)$ and for $y>1$ and any squarefree $d\ge1$,  we have
  \begin{equation*}
    |m_d(y)|\le g_0(d)\sqrt{2/y}
    +0.0144g_1(d)\frac{\1_{y\ge 10^{12}}}{\log y}
  \end{equation*}
  where the multiplicative functions $g_0$ and $g_1$ are defined on
  primes by:
  \begin{equation*}
    g_0(p)=
    \begin{cases}
      \sqrt{3/2}&\text{when $p=2$}\\
      \frac{\sqrt{p}}{\sqrt{p}-1}&\text{when $p\ge3$}
    \end{cases}
    \quad\text{and}\quad
    g_1(p)=
    \begin{cases}
      2.06&\text{when $p=2$}\\
      \frac{p^\xi}{p^\xi-1}&\text{when $p\ge3$}.
    \end{cases}
  \end{equation*}
\end{lem}
\noindent
This is a trivial modification of \cite[Proposition
5.15]{Helfgott*30}: when $d$ is even we remove only the coprimality to
$d/2$ and use directly the estimate for $m_2$ given in Lemma~\ref{m3}.
The value $2.06$ is an upper bound for $0.0296/0.0144$. The main
outcome of Lemma~\ref{m4} over \cite[Proposition
5.15]{Helfgott*30} is the improved value of $g_0(2)$.

\section{Auxiliaries}

Here is \cite[Lemma 3.2]{Ramare*13d}.
\begin{lem}
  \label{spe}
  Let $f$ be a C${}^1$ non-negative, non-increasing function over
  $[P,\infty[$, where $P\ge 3\,600\,000$ is a real number and such that $\lim_{t\rightarrow\infty}tf(t)=0$.
  We have
  \begin{equation*}
    \sum_{p\ge P} f(p)\log p
    \le (1+\epsilon) \int_P^\infty f(t) dt  +  \epsilon P f(P)  +  P f(P) / (5 \log^2 P)
  \end{equation*}
  with $\epsilon=1/914$. When we can only ensure $P\ge2$, then a similar
  inequality holds, simply replacing the last $1/5$ by a 4.
\end{lem}

\begin{lem}
  \label{Aux1}
  
  When $D\ge 0$, we have
  $\displaystyle
    \sum_{d\le D}\frac{\mu^2(d)\varphi(d)}{d}g_0(d)^2
    \le 2.07\, D$.
\end{lem}
\noindent
The maximum is reached at $D=42$.

\begin{lem}
  \label{Aux2}
  When $D\ge 0$, we have
  $\displaystyle
    \sum_{d\le
      D}\frac{\mu^2(d)\varphi(d)}{d}g_0(d)g_1(d)
    \le 1.60\, D$.
\end{lem}
\noindent
The maximum is reached at $D=7$.

\begin{lem}
  \label{Aux3}
  When $D\ge 0$, we have
  $\displaystyle
    \sum_{d\le D}\frac{\mu^2(d)\varphi(d)}{d}g_1(d)^2
    \le 1.57\,D$.
\end{lem}
\noindent
The maximum is reached at $D=3$.

\begin{proof}[Proof of Lemmas~\ref{Aux1},~\ref{Aux2} and~\ref{Aux3}]
  The three proofs are similar. We use $G$ for either the
  multiplicative function $g_0^2$, $g_0g_1$ or $g_1^2$.
  We readily check that
  \begin{equation*}
    \sum_{d\ge 1}
    \frac{\mu^2(d)\varphi(d)G(d)}{d^{1+s}}
    =\prod_{p\ge2}\biggl(
    1+\frac{(p-1)G(p)-p}{p^{1+s}}
    -\frac{(p-1)G(p)}{p^{1+2s}}
    \biggr)\zeta(s)
    =H(s)\zeta(s)
  \end{equation*}
  say. Notice that in the three cases, we have $(p-1)G(p)-p\ge0$.
  Thus, by adopting an obvious notation and using \cite[Lemma
  3.2]{Ramare*95} with $k_n=1/n$ and
  $g(d)={\mu^2(d)\varphi(d)G(d)}/d$ together with
  the second part of\cite[Lemma 2.1]{Ramare*14-1},
  we deduce that our sum $S$ satisfies
  \begin{equation*}
    S=
    H(1)D+\mathcal{O}^*(2.5\times\gamma \overline{H}(2/3)D^{2/3})
  \end{equation*}
  where
  \begin{equation*}
    H(1)=\prod_{p\ge2}\biggl(
    1+\frac{(p-1)G(p)-p}{p^2}
    -\frac{(p-1)G(p)}{p^{3}}
    \biggr)\le
    \begin{cases}
         2.0004&\text{when $G=g_0^2$,}\\
         1.34 &\text{when $G=g_0g_1$,}\\
         1.06 &\text{when $G=g_1^2$,}
    \end{cases}
  \end{equation*}
  and
  \begin{equation*}
    \overline{H}(2/3)=\prod_{p\ge2}\biggl(
    1+\frac{(p-1)G(p)-p}{p^{5/3}}
    +\frac{(p-1)G(p)(p)}{p^{7/3}}
    \biggr)\le
     \begin{cases}
         72.9&\text{when $G=g_0^2$,}\\
         23.4&\text{when $G=g_0g_1$,}\\
         9.20&\text{when $G=g_1^2$}.
    \end{cases}
  \end{equation*}
  Bounding above numerically $\overline{H}(2/3)$ requires some care:
  we use an Euler product for $p\le 10^8$ and the following Pari/GP-script
  which relies on Lemma~\ref{spe}:
\begin{verbatim}
{g0(p) = if(p == 2, return(sqrt(3/2)), return(1/(1-1/sqrt(p))));}

f0(p) = ((p-1)*g0(p)^2-p)/p^2 - (p-1)*g0(p)^2/p^3;
f1(p) = ((p-1)*g0(p)^2-p)/p^(5/3) + (p-1)*g0(p)^2/p^(7/3);

{val(boundP, myf) = 
   my(res = prodeuler(p = 2, boundP, 1.0 + myf(p)), eps = 1/914, aux);
   aux = (1+eps) * intnum(t = boundP, oo, myf(t)/log(t));
   aux += eps * boundP * myf(boundP) / log(boundP);
   aux += boundP * myf(boundP) / 5 /(log( boundP)^3);
   return(res * exp(aux));} 
\end{verbatim}
  We called \texttt{val(10000, f0)} and \texttt{val(10\^{}7, f1)}.
  We called \texttt{val(100000, f2)} and \texttt{val(10\^{}7, f3)}.
  We called \texttt{val(100000, f4)} and \texttt{val(10\^{}7, f5)}.
  This proves that, when $D\ge0$, we have
  \begin{equation*}
    \sum_{d\le D}\frac{\mu^2(d)\varphi(d)}{d}G(d)
    \le
    \begin{cases}
         2.0004 D + 106 D^{2/3}&\text{when $G=g_0^2$,}\\
         1.34 D +  33.8 D^{2/3}&\text{when $G=g_0g_1$,}\\
         1.06 D +  13.3 D^{2/3}&\text{when $G=g_1^2$.}
    \end{cases}
  \end{equation*}
  To complete, we called respectively
  \texttt{check(4*10\^{}9, af1atp, 2.0004, 106)},
   \texttt{check(10\^{}7, af2atp, 1.34, 33.8)} and
   \texttt{check(10\^{}6, af3atp, 1.06, 13.3)} of the same script.
\end{proof}

\section{A bound for the tail}

\begin{lem}
  \label{Le1}
  When $x\ge D\ge0$, we have
  \begin{equation*}
    \sum_{ d\le \min(D, x/10^{12})}
    \frac{\mu^2(d)\varphi(d)}{d^{3/2}\log(x/d)}g_0(d)g_1(d)
    \le 0.05\sqrt{D}.
  \end{equation*}
\end{lem}

\begin{proof}
    Set $y=\min(D, x/10^{12})$. We have
  \begin{equation*}
    \frac{d}{dt}\frac{1}{\sqrt{t}\log(x/t)}=
    -\frac{1}{2t^{3/2}\log (x/t)}+\frac{1}{t^{3/2}\log^2 (x/t)}
  \end{equation*}
  which is negative when $x/t\ge 10^{12}$.
  Hence, by Lemma~\ref{Aux2}, our sum $S$ satisfies
  \begin{align*}
    S
    &\le \int_1^y 1.60 \frac{2\log (x/t)-1}{2\sqrt{t}\log(x/t)^2}dt
    + \frac{1.60\sqrt{D}}{\log(x/D)}
    \\&\le 0.80\sqrt{x}\int_{x/y}^x  (2\log u-1)\frac{du}{u^{3/2}\log^2u}
    + 0.0497\sqrt{D}\le 0.05\sqrt{D}.
  \end{align*}
  The lemma follows readily.
\end{proof}

\begin{lem}
  \label{Le2}
  When $x\ge D\ge0$, we have
  \begin{equation*}
    \sum_{ d\le \min(D, x/10^{12})}
      \frac{\mu^2(d)\varphi(d)}{d^2\log(x/d)^2}g_1(d)^2
    \le 0.047.
  \end{equation*}
\end{lem}

\begin{proof}
  Set $y=\min(D, x/10^{12})$. We have
  \begin{equation*}
    \frac{d}{dt}\frac{1}{t\log(x/t)^2}=
    -\frac{1}{t^2\log (x/t)^2}+\frac{2}{t^2\log^2 (x/t)^3}
  \end{equation*}
  which is negative when $x/t\ge 10^{12}$.
  Hence, by Lemma~\ref{Aux3}, our sum $S$ satisfies
  \begin{align*}
    S
    &\le \int_1^y 1.57 \frac{\log (x/t)-2}{t\log^3(x/t)}dt
    + \frac{1.57}{\log(x/D)^2}
    \\&\le 1.57\int_{x/y}^x  (\log u-1)\frac{du}{u\log^3u}
    + 0.00152
    \le 0.047
  \end{align*}
  The lemma follows readily.
\end{proof}

\begin{lem}
  \label{Tail}
  When $x\ge D>0$, we have
  \begin{equation*}
    \sum_{ d\le D}\frac{\mu^2(d)\varphi(d)}{d^2}
    m_d(x/d)^2
    \le 4.14 \frac{D}{x}+0.00205.
  \end{equation*}
\end{lem}

\begin{proof}
  We readily find that
  \begin{align*}
    \sum_{ d\le D}\frac{\mu^2(d)\varphi(d)}{d^2}
    &\biggl(g_0(d)\sqrt{\frac{2d}{x}}
    +\1_{x/d\ge 10^{12}}g_1(d)\frac{0.0144}{\log (x/d)}\biggr)^2
    \\\le\, 
    &2.07\frac{2D}{x}
      +\frac{2\sqrt{2}\cdot 0.0144}{\sqrt{x}}
      \sum_{ d\le \min(D, x/10^{12})}
      \frac{\mu^2(d)\varphi(d)}{d^{3/2}\log(x/d)}g_0(d)g_1(d)
    \\
    &+0.0144^2
      \sum_{ d\le \min(D, x/10^{12})}
      \frac{\mu^2(d)\varphi(d)}{d^{2}\log(x/d)^2}g_1(d)^2.
  \end{align*}
  By Lemmas~\ref{Aux1},~\ref{Le1} and~\ref{Le2}, we get
  \begin{multline*}
    \sum_{ d\le D}\frac{\mu^2(d)\varphi(d)}{d^2}
    \biggl(g_0(d)\sqrt{\frac{2d}{x}}
    +\1_{x/d\ge 10^{12}}g_1(d)\frac{0.0144}{\log (x/d)}\biggr)^2
    \\\le
    4.14\frac{D}{x}
    +0.00204
    +0.00000975.
  \end{multline*}
\end{proof}

\section{Some refinements due to coprimality}

We shall use this lemma when bounding $|r_2^*(X;q)|$ below.
\begin{lem}
  \label{auxmajorstar2}
  We have
  \begin{equation*}
    \max_{M\ge1, K>0}\biggl|K\sum_{\substack{k\ge K \\
        (k,M)=1}}\frac{\mu(k)\varphi(k)}{k^3}\biggr|=1.
  \end{equation*}
\end{lem}

\begin{proof}
  Let us denote by $S$ our sum.
  \paragraph{When $K<1$:}
  We  notice that
  \begin{equation*}
    \sum_{\substack{k\ge 1 \\ (k,M)=1}}\frac{\mu(k)\varphi(k)}{k^3}
    =\prod_{p\neq M}\biggl(1-\frac{1}{p^2}+\frac{1}{p^3}\biggr)\le 1
  \end{equation*}
  which establishes the estimate $|S|\le 1/K$ when $K\le 1$. 
  \paragraph{When $1\le K<2$:}
  We find that
  \begin{equation*}
    S = \prod_{p\neq M}\biggl(1-\frac{1}{p^2}+\frac{1}{p^3}\biggr)-1.
  \end{equation*}
  Our sum is thus non-positive and its smallest value is
  \begin{equation*}
    \prod_{p\ge2}\biggl(1-\frac{1}{p^2}+\frac{1}{p^3}\biggr)-1\ge
    -0.252\ge -\frac{0.504}{K}\quad(\text{when $1<K\le 2$}).
  \end{equation*}
  \paragraph{When $2\le K<3$:}
  We find that
  \begin{equation*}
    S =
    \begin{cases}
      \prod_{p\neq M}\bigl(1-\frac{1}{p^2}+\frac{1}{p^3}\bigr)-1
      &\text{when $2|M$},\\
      \prod_{p\neq M}\bigl(1-\frac{1}{p^2}+\frac{1}{p^3}\bigr)-\frac78
      &\text{when $2\nmid M$}.
    \end{cases}
  \end{equation*}
  This implies in the first case that $S$ is non-positive, and that it
  bounded above by $1/8$ in the second case. Thus
   \begin{equation*}
    |S| \le 
    \begin{cases}
      0.145\le \frac{0.435}{K}
      &\text{when $2|M$},\\
      \max(0.127,1/8)\le \frac{0.381}{K}
      &\text{when $2\nmid M$}.
    \end{cases}
  \end{equation*}
  \paragraph{When $2|M$:}
  In that case, $k$ is odd and a comparison to an integral gives us
  \begin{equation}
    \label{useful}
    |S|\le \frac{1}{K^2}+\int_{(K-1)/2}^\infty\frac{dt}{(2t+1)^2}
    \le \frac{1}{K^2}+\frac{1}{2K}\le \frac{5}{6K}
  \end{equation}
  on assuming $K\ge3$, establishing our estimate in this case.
  \paragraph{When $3\le K<4$ and $2\nmid M$:}
  We find that
  \begin{equation*}
    S =
    \begin{cases}
      \prod_{p\neq M}\bigl(1-\frac{1}{p^2}+\frac{1}{p^3}\bigr)-\frac78
      &\text{when $3|M$},\\
      \prod_{p\neq M}\bigl(1-\frac{1}{p^2}+\frac{1}{p^3}\bigr)-\frac{173}{216}
      &\text{when $3\nmid M$}.
    \end{cases}
  \end{equation*}
  This implies that
   \begin{equation*}
    |S| \le 
    \begin{cases}
      \max(\frac{25}{27}-\frac78, 0.127)\le \frac{0.508}{K}
      &\text{when $3|M$},\\
      0.053\le \frac{0.212}{K}
      &\text{when $3\nmid M$}.
    \end{cases}
  \end{equation*}
  This proves our estimate in this case.
  \paragraph{When $2\nmid M$:} 
  We may write
  \begin{align*}
    S
    &=\sum_{\substack{k\ge K\\ (k, M)=1\\ (k,2)=1}}\frac{\mu(k)\varphi(k)}{k^3}
    +\sum_{\substack{k\ge 2K\\ (k, M)=1\\ 2|k}}\frac{\mu(k)\varphi(k)}{k^3}
    +\sum_{\substack{2K > k\ge K\\ (k, M)=1\\ 2|k}}\frac{\mu(k)\varphi(k)}{k^3}
    \\&=\frac{7}{8}\sum_{\substack{k\ge K\\ (k, M)=1\\ (k,2)=1}}\frac{\mu(k)\varphi(k)}{k^3}
    -\frac18 \sum_{\substack{K > k\ge K/2\\ (k, M)=1\\ (k,2)=1}}\frac{\mu(k)\varphi(k)}{k^3}
  \end{align*}
  We use Eq.~\eqref{useful} to infer that
  \begin{equation*}
    |S|\le \frac{7}{8} \biggl(\frac{1}{K^2}+\frac{1}{2K}\biggr)
    + \frac{1}{8} \biggl(\frac{1}{(K/2)^2}+\frac{1}{2(K/2)}\biggr).
  \end{equation*}
  We readily check that this quantity is $ < 1/ K$ when $K\ge 4$,
  therefore concluding this proof. 
\end{proof}

\section{On a family of functions: initial step}

In this section and the next one, we study the family
\begin{equation}
  \label{defGstarq}
  G^*_q(X)=\sum_{\substack{d\le X\\
      (d,q)=1}}\frac{\mu^2(d)\varphi(d)}{d^2},
  \quad
  G^*(X)=G^*_1(X).
\end{equation}
It will transpire from the proof to come that the study of $G^*$ is of
special importance for the general case. We devote this section to
the precise modification of initial case that will be required. The
approach is different: in the next section, we
shall compare the function $\1_{(d,q)=1}\mu^2(d)\varphi(d)/d^2$ to the
function $1/d$ (see Lemma~\ref{convol}) while here, we compare the function
$\mu^2(d)\varphi(d)/d^2$ to $\mu^2(d)/d$, as in \cite{Ramare*18-9}. Such a comparison is the
subject of our first lemma.
\begin{lem}
  \label{convol0}
  We have
  \begin{equation*}
    \mu^2(d)\frac{\varphi(d)}{d}
    =
    \sum_{\ell m=d}\mu^2(\ell)g(m)
  \end{equation*}
  where $g$ is the multiplicative function defined by
  $g(p^k)=(-1)^k/p$ for every positive integer $k$ and every prime~$p$.
\end{lem}

\begin{proof}
  We simply compare the $p$-factors of the corresponding Dirichlet
  series, and check that
  \begin{equation*}
    \frac{1+\frac{p-1}{p^{1+s}}}{1+\frac{1}{p^s}}
    =
    1-\frac{1}{p^s}\frac{1}{1+\frac{1}{p^s}}
    =1+\sum_{k\ge1}\frac{g(p^k)}{p^{ks}}
  \end{equation*}
  from which the lemma follows.
\end{proof}

\begin{lem}
  \label{MoebiusSquare}
  We have
  \begin{equation*}
    \forall X\ge X_0,\quad
    \sum_{d\le X}\mu^2(d)=\frac{6}{\pi^2}X+\Ocal^*\bigl(c(X_0)\sqrt{X}\bigr)
  \end{equation*}
  where
  \begin{tabular}{|c|c|c|c|c|c|}
    \hline
    $X_0$& 0&8& 1664&82005&438653\cr
    \hline
    $c(X_0)$& 1&1/2&0.1333&0.036438&0.02767\cr
    \hline
  \end{tabular}.
\end{lem}

\begin{lem}
  \label{init}
  We have
  \begin{equation*}
    \max_{X\ge 0}\frac{1}{\sqrt{X}}
    \biggl(\sum_{d\le X}\mu^2(d)\frac{\varphi(d)}{d}-AX\biggr)=
    1-A\le 0.572
  \end{equation*}
  where
  \begin{equation*}
    A = \prod_{p\ge2}\biggl(1-\frac{2}{p^2}+\frac{1}{p^3}\biggr)
    =0.428257\cdots.
  \end{equation*}
\end{lem}

\begin{proof}
  By Lemma~\ref{convol0} and when $X\ge M_0X_0$, we may write
  \begin{align*}
    \sum_{d\le X}\mu^2(d)\frac{\varphi(d)}{d}
    &=\sum_{m \ge 1} g(m)\sum_{\ell\le X/m}\mu^2(\ell)
    \\&=\sum_{m \le M_0}
    g(m)\biggl(\frac{6}{\pi^2}\frac{X}{m}+\Ocal^*\biggl(c(X_0)\sqrt{\frac{X}{m}}\biggr)\biggr)
    \\&\qquad+\sum_{m >M_0 }
    g(m)\biggl(\frac{6}{\pi^2}\frac{X}{m}+\Ocal^*\biggl(\sqrt{\frac{X}{m}}\biggr)\biggr)
    \\ &=AX
         +
         \Ocal^*\biggl(\sqrt{X}\sum_{m\ge1}\frac{|g(m)|}{\sqrt{m}}
         \times
         \begin{cases}
           c(X_0)&\text{when $m\le M_0$},\\
           1&\text{when $m> M_0$}.
         \end{cases}
             \biggr)
  \end{align*}
  We have therefore reached the fundamental formula, valid for $X\ge X_0M_0$:
  \begin{equation*}
    \sum_{d\le X}\mu^2(d)\frac{\varphi(d)}{d}
    =A X
    +\Ocal^*\biggl(\biggl(\sum_{m\ge1}\frac{|g(m)|}{\sqrt{m}}
    -(1-c(X_0))\sum_{m\le M_0}\frac{|g(m)|}{\sqrt{m}}\biggr)\sqrt{X}\biggr).
  \end{equation*}
  We numerically check that
  \begin{equation*}
    \forall X\le 4\cdot 10^9,\quad
    \sum_{d\le X}\mu^2(d)\frac{\varphi(d)}{d}\le
    AX+(1-A)\sqrt{X}.
  \end{equation*}
  See \texttt{AMoebiusSum-r1-01.gp/getmaxlocr0()}. Therefore, we may
  take $M_0=9118$. The script
  \texttt{AMoebiusSum-r1-01.gp/getmaxasympr0()} concludes the proof.
\end{proof}

\section{On a family of functions: first batch}

Let us mention that a related study appears in \cite[Lemma
4.5]{Alterman*22} by S. Zuniga Alterman.
Here is the final result of this part.
We proceed as for the proof of \cite[Theorem 1.1]{Ramare*14-1} which
relies on \cite[Theorem 1.4]{Ramare*14-1}. Let us start by a more
precise form though more special form of this lemma.
\begin{lem}
  \label{keyb}
  Let $(g(m))_{m\ge1}$ be a sequence of complex numbers such that both series
  $\sum_{m\ge1} g(m)/m$ and $\sum_{m\ge1} g(m)(\log m)/m$ converge. We define
  $G^\sharp(x)=\sum_{m> x} g(m)/m$ and assume that
  $\int_1^\infty |G^\sharp(t)|dt/t$ converges.
  We then have
  \begin{multline*}
    \sum_{n\le D}\frac{(g\star\1)(n)}{n}
    =\sum_{m\ge1}\frac{g(m)}{m}\Bigl(\log\frac{D}{m}+\gamma\Bigr)
    +\int_{e^\gamma D}^\infty G^\sharp(t)\frac{dt}{t}
    \\
    +\Ocal^*\Bigl(\frac{1}{D}\int_{1}^{e^\gamma}\sum_{m\le u D}|g(m)|\frac{du}{u}\Bigr).
  \end{multline*}
\end{lem}

\begin{proof}
  We established in \cite[Eq.~(2.2)]{Ramare*14-1} the following formula:
  \hspace{0.1cm}\par\noindent
  \fbox{\begin{minipage}{0.9\linewidth}\begin{multline}\label{iden1}
        \sum_{n\le D}\frac{(g\star \1)(n)}{n}
        =
        \sum_{m\ge1} \frac{g(m)}{m} \Bigl(\log\frac{D}{m}+\gamma\Bigr)
        +\int_{\eta D}^\infty  G^\sharp(t)\frac{dt}{t}
        \\
        -(\gamma-\log\eta) G^\sharp(\eta D)
        +
        \sum_{m\le \eta D} \frac{g(m)}{m} R\left(\frac{D}{m}\right).
      \end{multline}
    \end{minipage}
  }
  
  \smallskip
  \noindent
  where $\eta\ge1$ and $R(t)=\sum_{n\le t}1/n-\log t-\gamma$. We also
  showed that, for any positive $t$, we have $|R(t)|\le \gamma$. Let us
  specialize $\eta=e^\gamma$. We see that
  $R\left(\frac{D}{m}\right)=-\int_{D/m}^{e^\gamma}du/u$, so that
  \begin{equation*}
    \sum_{m\le \eta D} \frac{|g(m)|}{m}
    \biggl|R\left(\frac{D}{m}\right)\biggr|
    \le \int_1^{e^\gamma}\sum_{m\le D}\frac{|g(m)|}{m}\frac{du}{u}
    + \int_1^{e^\gamma}\sum_{D<m\le uD}\frac{|g(m)|}{m}\frac{du}{u}
  \end{equation*}
  from which our lemma readily follows.
\end{proof}

\begin{lem}
  \label{GetGstarq}
  For every $X>0$, we have
  \begin{equation*}
    G^*_q(X)=
    A\prod_{p|q}\frac{p^2}{p^2+p-1}\bigl(\log X+c_q\bigr)
    +\Ocal^*(4.73 j_1^*(q)/\sqrt{X})
  \end{equation*}
  where (as in Lemma~\ref{init})
  \begin{equation*}
    A = \prod_{p\ge2}\biggl(1-\frac{2}{p^2}+\frac{1}{p^3}\biggr)
    =0.428257\cdots,
    \quad
     j_1^*(q)=\prod_{p|q}\frac{p^{3/2}+p}{p^{3/2}+1},
  \end{equation*}
  then
  \begin{equation*}
    c_q=\gamma+\sum_{p|q}\frac{(p-1)\log p}{p^2+p-1}
    +
    \sum_{p\ge 2}\frac{(3p-2)\log p}{(p-1)(p^2+p-1)},
  \end{equation*}
  Morevover
  \begin{multline*}
    G^*_q(X)-G_q^*(Y)=
    A\prod_{p|q}\frac{p^2}{p^2+p-1}\log \frac{X}{Y}
    \\+\Ocal^*\biggl(2.18 j_1^*(q)
    \biggl(
    \frac{2 (e^{\gamma/2}-1)}{\sqrt{X}}
    +\frac{2 (e^{\gamma/2}-1)}{\sqrt{Y}}
    +\frac{2}{\sqrt{e^\gamma Y}}
    -\frac{2}{\sqrt{e^\gamma X}}
      \biggr)\biggr).
  \end{multline*}
\end{lem}

\begin{lem}
  \label{convol}
  We have
  \begin{equation*}
    \1_{(d,q)=1}\frac{\mu^2(d)\varphi(d)}{d}
    =
    \sum_{\substack{k^2\ell r|d \\ r|q\\ (k\ell,q)=1\\ (k,\ell)=1}}\frac{\mu(rk\ell)\varphi(k)}{k\ell}.
  \end{equation*}
\end{lem}
\noindent
This is the counterpart of \cite[Lemma 4.1]{Ramare*14-1}.
\begin{proof}
  We find that
  \begin{equation*}
    D_q(s)=\sum_{\substack{d\ge 1\\
        (d,q)=1}}\frac{\mu^2(d)\varphi(d)}{d^{1+s}}
    =\prod_{\substack{p\ge 1\\ (p,q)=1}}\biggl(
    1+\frac{p-1}{p^{1+s}}
    \biggr)
  \end{equation*}
  which we decompose in
  \begin{equation}
    \label{decDq}
    D_q(s)=\zeta(s)\prod_{p|q}\biggl(1-\frac{1}{p^s}\biggr)
    \prod_{\substack{p\ge 1\\ (p,q)=1}}\biggl(
    1-\frac{1}{p^{1+s}}-\frac{p-1}{p^{1+2s}}
    \biggr).
  \end{equation}
  The identity then follows either by using the unitarian convolution
  as in \cite[Theorem 3.5]{Ramare*22-0} or by simply checking that
  both left and right hand side of are multiplicative functions, and
  that they coincide on prime powers.
\end{proof}
We define in this section
\begin{equation}
  \label{defr2star}
  r_2^*(X;q)=\sum_{\substack{k^2\ell r>X\\ r|q\\ (k\ell,q)=(k,\ell)=1}}
  \frac{\mu(rk\ell)\varphi(k)}{rk^3\ell^2}
\end{equation}
as well as
\begin{equation}
  \label{defr1star}
  r_1^*(X;q)
  =\sum_{\substack{k^2\ell r\le X\\ r|q\\ (k\ell,q)=(k,\ell)=1}}
  \frac{\mu^2(rk\ell)\varphi(k)}{k\ell}.
\end{equation}
Let us first majorize $r_1^*(X;q)$.
\begin{lem}
  \label{majorstar1}
  The function $j_1^*$ being defined in Lemma~\ref{GetGstarq}, we have
  \begin{equation*}
    \forall X\ge 0,\quad r_1^*(X;q)\le 2.18\sqrt{X}j_1^*(q).
  \end{equation*}
  We also have  
  \begin{equation*}
    r_1^*(X;q)\le 0.931\sqrt{X}j_1^*(q)
    +1.96 X^{1/4} j_5^*(q).
  \end{equation*}
  where
  \begin{equation}
    \label{eq:5}
    j_5^*(q)=\prod_{p|q}\frac{p^{5/4}+p}{p^{5/4}+1}.
  \end{equation}
\end{lem}
\noindent
This is the counterpart of \cite[Lemma 6.1]{Ramare*14-1}.
\begin{proof}
  We take advantage of the variable $k$ by writing
  \begin{equation*}
    r_1^*(X;q)
    \le
      \sum_{\substack{\ell r\le X\\ r|q\\ (\ell,q)=1}}
    \frac{\mu^2(r\ell)}{\ell}\sqrt{\frac{X}{\ell r}}
    =\sqrt{X}\prod_{p|q}\frac{p^{3/2}+p}{p^{3/2}+1}
    \prod_{p\ge2}\biggl(1+\frac{1}{p^{3/2}}\biggr).
  \end{equation*}
  We finally notice that
  $\prod_{p\ge2}(1+p^{-3/2})=\zeta(3/2)/\zeta(3)$. The first part of
  the lemma follows readily. The proof we have just followed used the
  upper bound
  \begin{equation*}
    \sum_{\substack{k\le K\\ (k,q\ell)=1}}
    \mu^2(k)\frac{\varphi(k)}{k}\le K. 
  \end{equation*}
  As it turns out, the quantity to be majorized can be handled by
  Lemma~\ref{GetGstarq}! Such a recursive treatment leads to constants
  that are too big for us. But we may still forget of the coprimality
  condition and appeal to Lemma~\ref{init}. This gives us
  \begin{align*}
    r_1^*(X;q)
    &\le
      A\sum_{\substack{\ell r\le X\\ r|q\\ (\ell,q)=1}}
    \frac{\mu^2(r\ell)}{\ell}\sqrt{\frac{X}{\ell r}}
    +(1-A)
    \sum_{\substack{\ell r\le X\\ r|q\\ (\ell,q)=1}}
    \frac{\mu^2(r\ell)}{\ell}\biggl(\frac{X}{\ell r}\biggr)^{1/4}
    \\&\le
    A\frac{\zeta(3/2)}{\zeta(3)}j_1^*(q)\sqrt{X}
    + (1-A)\frac{\zeta(5/4)}{\zeta(5/2)}
    \prod_{p|q}\frac{p^{5/4}+p}{p^{5/4}+1}.
  \end{align*}
  A numerical application concludes the proof.
\end{proof}

\begin{lem}
  \label{majorstar2}
  The function $j_1^*$ being defined in Lemma~\ref{GetGstarq}, we have
  \begin{equation*}
    |r^*_2(X;q)|\le \frac{2.18}{\sqrt{X}}j_1^*(q).
  \end{equation*}
\end{lem}
\noindent
This is the counterpart of \cite[Lemma 6.2]{Ramare*14-1}.
\begin{proof}
  We again take advantage of the variable $k$ and use Lemma~\ref{auxmajorstar2}.
  \begin{align*}
    |r^*_2(X;q)|
    &\le
    \sum_{\substack{\ell\ge 1, r|q\\ (\ell,q)=1}}
    \frac{\mu^2(r\ell)}{r\ell^2}
    \biggl|\sum_{\substack{k^2 >X/(\ell r)\\ (k,q\ell)=1}}
    \frac{\mu(k)\varphi(k)}{k^3}\biggr|
    \\&\le \frac{1}{\sqrt{X}}
    \sum_{\substack{\ell\ge 1, r|q\\ (\ell,q)=1}}
    \frac{\mu^2(r)\mu^2(\ell)}{\sqrt{r}\ell^{3/2}}
  \end{align*}
  where we recognize the quantities that appeared in the proof of
  Lemma~\ref{majorstar1}. The lemma then follows swiftly.
\end{proof}

\begin{proof}[Proof of Lemma~\ref{GetGstarq}]
  We employ \cite[Theorem 1.4]{Ramare*14-1} with the function $g$
  being
  \begin{equation*}
    g(m)=\sum_{\substack{k^2\ell r=m \\ r|q\\ (k\ell,q)=1\\
        (k,\ell)=1}}
    \frac{\mu(rk\ell)\varphi(k)}{k\ell}.
  \end{equation*}
  This is a consequence of Lemma~\ref{convol}. 
  Then $r_2^*(X;q)$ is
  the $G^\sharp(X)$ of \cite[Theorem 1.4]{Ramare*14-1} while
  $r_1^*(X;q)$ is $\sum_{m\le X}|g(m)|$. As a consequence, and with
  $\eta=e^\gamma$, we find that, for $X\ge1$, we have
  \begin{equation}
    \label{raw}
    G^*_q(X)=
    \sum_{m\ge 1}\frac{g(m)}{m}\Bigl(\log\frac{X}{m}+\gamma\Bigr)
    +\int_{\eta X}r_2^*(t;q)\frac{dt}{t}
    +\Ocal^*\biggl(\int_1^{e^\gamma} \frac{r_1^*(u X;q)du}{uX}\biggr).
  \end{equation}
  Lemmas~\ref{majorstar1} and~\ref{majorstar2} gives
  \begin{equation*}
    G^*_q(X)=
    \sum_{m\ge 1}\frac{g(m)}{m}\Bigl(\log\frac{X}{m}+\gamma\Bigr)
    +\Ocal^*\biggl(2.18(2 (e^{\gamma/2}-1)+2e^{-\gamma/2}) \frac{j_1^*(q)}{\sqrt{X}}\biggr).
  \end{equation*}
  We identify the main term by using the Dirichlet series
  $H_q(s)=D_q(s)/\zeta(s)$ (see Eq.~\eqref{decDq}) of~$g$, and get
  \begin{equation*}
    G^*_q(X)=
    H_q(1)\biggl(\log X+\frac{H_q'(1)}{H_q(1)}+\gamma\Bigr)
    +\Ocal^*(4.73 j_1^*(q)/\sqrt{X}).
  \end{equation*}
  We find that
  \begin{equation*}
    H_q(1)
    =
    \prod_{p|q}\frac{p^2}{p^2+p-1}
    \prod_{p\ge2}\biggl(1-\frac{2}{p^2}+\frac{1}{p^3}\biggr).
  \end{equation*}
  We get to $H'_q/H_q$ by using the logarithmic derivatives and find
  that
  \begin{align*}
    \frac{H_q'(1)}{H_q(1)}
    &=\sum_{p|q}\frac{\log p}{p-1}
    +
    \sum_{p\nmid q}\frac{(3p-2)\log p}{(p-1)(p^2+p-1)}
    \\&=\sum_{p|q}\frac{(p-1)\log p}{p^2+p-1}
    +
    \sum_{p\ge 2}\frac{(3p-2)\log p}{(p-1)(p^2+p-1)}.
  \end{align*}
  \paragraph{On the difference $G^*_q(X)-G^*_q(Y)$:}
  To get the more precise evaluation of $G^*_q(X)-G^*_q(Y)$, we go
  back to Eq.~\eqref{raw} to save on the factor involving $r_2^*$.
  %
  The
  proof is then rapidly completed.
\end{proof}

\section{On a family of functions: second batch}

\begin{lem}
  \label{major1starter}
  When $q$ has all its prime factors below 30, we have
  \begin{equation*}
    \forall X\ge 0,\quad r_1^*(X; q)\le 1.17\sqrt{X}j_1^*(q).
  \end{equation*}
\end{lem}

\begin{proof}
  We first prove by hard computations,  that
  when $q$ has all its prime factors below 30, we have
  \begin{equation*}
    \forall X\le 10^6,\quad r_1^*(X; q)\le 1.17\sqrt{X}j_1^*(q).
  \end{equation*}
  See the Pari/GP script
  \texttt{AMoebiusSum-r1-01.gp/getmaxr1()}.  
  For $X$ large, we use the second bound provided by
  Lemma~\ref{majorstar1} and get
  \begin{equation*}
    \frac{r_1^*(X;q)}{\sqrt{X}j_1^*(q)}\le
    0.931
    +1.96\frac{j_5^*(q)}{X^{1/4}j_1^*(q)}.
  \end{equation*}
  The lemma follows readily.
\end{proof}

\section{Direct computations}

\begin{lem}
  When $422\le X\le 11\,000\,000$, we have 
  \begin{equation*}
    0\le \sum_{\substack{d_1,d_2\le
        X}}\frac{\mu(d_1)\mu(d_2)}{[d_1,d_2]}
    \le 0.445.
  \end{equation*}
  On $[6, 10\,040\,000]$, this sum is bounded above by $0.528$.
  On $[2, 10\,040\,000]$, this sum is bounded above by
  $19/30=0.633\cdots$.
  When $X\ge 1000$, this sum remained $\ge 0.437$.
\end{lem}
A value larger than $0.44455$ is reached around $D=1321$.
Let us set
\begin{equation}
  \label{eq:2}
  \Sigma(X)=\sum_{\substack{d_1,d_2\le
        X}}\frac{\mu(d_1)\mu(d_2)}{[d_1,d_2]}.
\end{equation}
When $X$ is an integer, say $d$, we find that
\begin{equation*}
  \Sigma(d)-\Sigma(d-1)
  =\frac{\mu^2(d)}{d}
  +2\mu(d)\sum_{d' < d}\frac{\mu(d')}{[d,d']}.
\end{equation*}
This yields a formula to find the maximum of $\Sigma(d)$ over some
range, but each step is costly. We continue with
\begin{align*}
  \Sigma(d)-\Sigma(d-1)
  &=\frac{\mu^2(d)}{d}
  +2\frac{\mu(d)}{d}\sum_{d' < d}(d,d')\frac{\mu(d')}{d'}.
  \\&=\frac{\mu^2(d)}{d}
    +2\frac{\mu(d)}{d}\sum_{\delta |d}\mu(\delta)
    \sum_{\substack{d' < d/\delta \\ (d',d)=1}}\frac{\mu(d')}{d'}.
\end{align*}
We finally use the Landau formula (see for instance \cite[(5.73)]{Helfgott*30}):
\begin{equation*}
  \sum_{\substack{d' < d/\delta \\ (d',d)=1}}\frac{\mu(d')}{d'}
  =\sum_{\ell|d^\infty}\frac{1}{\ell}\sum_{\substack{d' < d/(\ell \delta)}}\frac{\mu(d')}{d'}.
\end{equation*}
Therefore
\begin{equation*}
  \Sigma(d)-\Sigma(d-1)
  =
  \frac{\mu^2(d)}{d}
  +2\frac{\mu(d)}{d}\sum_{\substack{\delta |d\\ \ell|d^\infty}}\frac{\mu(\delta)}{\ell}
    m((d-1)/(\delta\ell)).
\end{equation*}
Let us join $\delta$ and $\ell$ in $k=\delta \ell$. We
have $k|d^\infty$ and 
\begin{equation*}
  \sum_{\delta \ell=k}\frac{\mu(\delta)}{\ell}=
  \frac{1}{k}\prod_{p|k}(1-p)=\frac{(-1)^{\Omega(k)}\varphi(k)}{k^2}.
\end{equation*}
Here is thus the identity we use:
\begin{equation}
  \label{eq:3}
  \Sigma(d)-\Sigma(d-1)
  =
  \frac{\mu^2(d)}{d}
  +2\frac{\mu(d)}{d}\sum_{\substack{k|d^\infty}}\frac{(-1)^{\Omega(k)}\varphi(k)}{k^2}
    m\Bigl(\frac{d-1}{k}\Bigr).
\end{equation}
See script \texttt{AMoebiusSum-02.gp/DITgreat()}.
This entails to precompute all the values $m(t)$ for $t$ upto where we
want to compute and this is a large array. So we decided to only store
the values $m_{M_0}(t)$ for $M_0=6$. We stored in fact
\begin{equation}
  \label{eq:1}
  m(t;u,M_0)=\sum_{\substack{n\le t\\ n\equiv u[M_0]}}\mu(n)/n
\end{equation}
for $u$  covering a reduced congruence system modulo~$M_0$, i.e. in
practice all
$u\in\{1,\cdots, M_0\}$ that are coprime to $M_0$. This reduced
sizeably the amount of storage while accessing to the values $m(t)$ takes more
time in a classical time/space bargain.
See script \texttt{AMoebiusSum-02.gp/DITb()}.

\section{Proof of Theorem~\ref{Main}}

\begin{proof}
  We readily find that
  \begin{align*}
    S
    &=\sum_{d\le x}\frac{\mu^2(d)\varphi(d)}{d^2}
    \biggl(\sum_{\substack{n\le x/d\\ (n,d)=1}}\frac{\mu(n)}{n}\biggr)^2
    \\&= \sum_{D< d\le x}\frac{\mu^2(d)\varphi(d)}{d^2}m_d(x/d)^2
    +\sum_{ d\le D}\frac{\mu^2(d)\varphi(d)}{d^2}m_d(x/d)^2.
  \end{align*}
  The second sum is handled in Lemma~\ref{Tail}. Let us call $S_0$ the
  first one.
  We set $\Delta(j)=\prod_{p\le j}p$ and write
  \begin{align*}
    S_0
    &=
    \sum_{j\le x/D}\sum_{\max(x/D,x/(j+1))< d\le x/j}\frac{\mu^2(d)\varphi(d)}{d^2}m_d(j)^2
    \\&=
    \sum_{j\le x/D}\sum_{\delta|\Delta(j)}
    m_\delta(j)^2
    \sum_{\substack{\max(x/D,x/(j+1))< d\le x/j\\ (d,\Delta(j))=\delta}}\frac{\mu^2(d)\varphi(d)}{d^2}
    \\&\le
    \sum_{j\le x/D}\sum_{\delta|\Delta(j)}
    \frac{\mu^2(\delta)\varphi(\delta)}{\delta^2}m_\delta(j)^2
    \sum_{\substack{\frac{x}{(j+1)\delta}< d\le \frac{x}{j\delta}
    \\ (d,\Delta(j))=1}}\frac{\mu^2(d)\varphi(d)}{d^2}.
  \end{align*}
  A direct usage of Lemma~\ref{GetGstarq} gives us
  \begin{multline*}
    S_0
    =
    A\sum_{j\le x/D}
    \prod_{p|\Delta(j)}\frac{p^2}{p^2+p-1}\log\frac{j+1}{j}
    \sum_{\delta|\Delta(j)}
    \frac{\mu^2(\delta)\varphi(\delta)}{\delta^2}m_\delta(j)^2
    \\+\Ocal^*\biggl(
    \frac{1}{\sqrt{x}}
    \sum_{j\le x/D}j_1^*(\Delta(j))
    \biggl(
    2\times2.18 (e^{\gamma/2}-1)(\sqrt{j+1}+\sqrt{j})
    \\+\frac{2\times 2.18 e^{-\gamma/2}}{\sqrt{j+1}+\sqrt{j}}
      \biggr)
      \sum_{\delta|\Delta(j)}
    \frac{\mu^2(\delta)\varphi(\delta)}{\delta^2}m_\delta(j)^2\sqrt{\delta}
    \biggr)  .
  \end{multline*}
  This expression may be improved in two manners: when $\delta$ has
  all its prime factors below 30, we may rely on
  Lemma~\ref{major1starter} rather than on Lemma~\ref{majorstar1},
  therefore changing the first 2.18 (in factor of
  $\sqrt{j+1}+\sqrt{j}$) by 1.17. The second improvement consists in
  localizing $X$ in an interval $[Y,2Y]$. When $j\delta > 2Y$, no
  contribution is to be incorporated.
  
  By selecting $x/D=22.99$ and assuming $x\ge 11\,000\,000$, we reach $S\le
  0.679$.
  See Pari/GP script
  \texttt{AMoebiusSumMT.gp/DoIt()}.
  \begin{itemize}
  \item When $x\ge  10^{9}$, we reach $S\le
    0.574$ (on taking $x/D=38.99$).
  \item When $x\ge 3\cdot 10^{10}$, we reach $S\le
    0.536$ (on taking $x/D=55.99$).
  \item When $x\ge 2.4\cdot 10^{12}$, we reach $S\le
    0.504$ (on taking $x/D=75.99$).
  \end{itemize}
\end{proof}


\end{document}